\documentclass[12pt,a4paper]{amsart}
\usepackage{amsmath, amssymb, amsfonts,enumerate}
\usepackage[all]{xy}
\usepackage{amscd}

\newtheorem{theorem}{Theorem}[section]

\newtheorem{corollary}[theorem]{Corollary}
\newtheorem{proposition}[theorem]{Proposition}

\newtheorem*{question*}{Question}

\theoremstyle{definition}

\newtheorem*{definition*}{Definition}

\theoremstyle{remark}
\newtheorem*{remark}{Remark}

\numberwithin{equation}{section}

\newcommand {\K}{\mathbb{K}} %% a field
\newcommand {\Z}{\mathbb{Z}} %% integers
\newcommand {\R}{\mathbb{R}} %% reals
\newcommand {\C}{\mathbb{C}} %% complex

\DeclareMathOperator{\GL}{GL}

\DeclareMathOperator{\Sym}{Sym}

\begin{document}
\title[On a characterization of locally finite groups]{On a characterization of locally finite groups in terms of linear cellular automata}
%\date{\today}

\author{Tullio Ceccherini-Silberstein}
\address{Dipartimento di Ingegneria, Universit\`a del Sannio, C.so
Garibaldi 107, 82100 Benevento, Italy}
\email{tceccher@mat.uniroma1.it}
\author{Michel Coornaert}
\address{Institut de Recherche Math\'ematique Avanc\'ee,
                                            Universit\'e  de Strasbourg et CNRS,
                                                 7 rue Ren\'e-Descartes,
                                               67084 Strasbourg Cedex, France  }
\email{coornaert@math.u-strasbg.fr}
% \date{\today, Preliminary version}
\begin{abstract}
We prove that a group $G$ is locally finite if and only if every surjective real (or complex) linear cellular automaton with finite-dimensional alphabet over $G$ is injective.
 \end{abstract}

\maketitle

% SECTION 1
\section{Introduction}

A group is called \emph{locally finite} if all its finitely generated subgroups are finite.
A group is called \emph{periodic} if all its elements have finite order.
In a locally finite group, all cyclic groups are finite. Therefore every locally finite group is periodic.
Conversely, it is well known that every periodic solvable group is locally finite (this follows from an immediate induction on the derived length of the group). 
 However, there exist periodic groups which are not locally finite.
In fact, the question asking whether every periodic group is locally finite is a way to formulate the \emph{general Burnside problem}, which was raised by W. Burnside in 1902  and answered negatively by E.S. Golod and I.R. Shafarevich in 1964. A famous example of a group which is periodic but not locally finite is provided by  the \emph{Grigorchuk group}
which is a finitely generated infinite amenable $2$-group discovered by R.I. Grigorchuk \cite{grigorchuk} in 1980. 
\par
The class of locally finite groups is closed under taking subgroups, taking quotients, taking extensions, and taking direct sums.
Many results on finite groups, such as the theory of Sylow subgroups (see e.g. \cite{kegel} and the references therein), have been extended to locally finite groups.
 Of course every finite group is locally finite but there are also interesting families of infinite locally finite groups arising naturally   in various branches of mathematics. 
 For example, if $X$ is an infinite set, then the group $\Sym_0(X)$ consisting of all  finitely-supported permutations of $X$ is an infinite locally finite group.
Another important family of infinite locally finite groups is given by the linear groups $\GL_n(\K)$,
where $\K$ is the algebraic closure of a finite field.
\par
The goal of this paper is to present a characterization of locally finite groups in terms of linear cellular automata.
Before stating our main results, let us first recall some basic definitions.
\par
Let $G$ be a group and let $A$ be a set.
The set $A^G = \prod_{g \in G} A$, which we regard as the set consisting of all maps 
$x \colon G \to A$, is called the set of \emph{configurations} over the group $G$ and the \emph{alphabet} $A$. 
We make $G$ act on  $A^G$  by setting
$gx(h) = x(g^{-1}h)$ for all $x \in A^G$ and $g,h \in G$.
A \emph{cellular automaton} over the group $G$ and the alphabet $A$ is a map 
$\tau \colon A^G \to A^G$ satisfying the following condition:
there exist a finite subset $M \subset G$ 
and a map $\mu \colon  A^M \to A$ such that 
\begin{equation*} 
  \tau(x)(g) = \mu((g^{-1}x)\vert_M) \quad \text{for all } x \in A^G \text{ and } g \in G,
\end{equation*}
where $(\cdot)\vert_M \colon A^G \to A^M$ denotes the restriction map.
  Such a set $M$ is called a \emph{memory set}
    and $\mu$ is called a \emph{local defining map}
   for the cellular automaton $\tau$.
\par
 Suppose now that $V$ is a vector space over  a field $\K$.
A cellular automaton $\tau \colon V^G \to V^G$ is called a \emph{linear cellular automaton} if $\tau$ is $\K$-linear with respect to the product vector space structure on $V^G$.
 Our characterization of locally finite groups  is based on the relationship between injectivity and surjectivity for linear cellular automata.
The study of this relationship, and more especially of the domains of validity of implications of the type ``injective $\Rightarrow$ surjective" and ``surjective $\Rightarrow$ injective", is a recurrent theme in the literature on cellular automata.
For example, the classical 
Moore-Myhill Garden of Eden theorem (cf. \cite{moore}, \cite{myhill}) asserts that a cellular automaton with finite alphabet over the group $\Z^2$ is surjective if and only if it is pre-injective.
As the terminology indicates, pre-injectivity is a weak form of injectivity.
It means that any two configurations with the same image must be equal if they coincide outside a finite subset of the group
 (this is  equivalent to the absence of mutually erasable patterns for the cellular automaton, cf. the appendix in \cite{gromov}). 
 The Moore-Myhill theorem was extended to amenable groups in \cite{ceccherini}.
There is also a version of the Garden of Eden theorem for linear cellular automata with 
finite-dimensional alphabet over amenable groups
in \cite{Eden-linear}.
The question asking whether every injective cellular automaton with finite alphabet is surjective 
remains open and is known as the \emph{Gottschalk conjecture}.
Recently, this conjecture was answered affirmatively for a very large class of groups, namely \emph{sofic groups}, by M. Gromov \cite{gromov} and B. Weiss \cite{weiss}. The class of sofic groups was  introduced by Gromov in \cite{gromov}. It includes in particular all residually amenable groups, and therefore all amenable groups and all residually finite groups 
(actually, there is no known example of a non-sofic group up to now).  
 When $V$ is a finite-dimensional vector space and $G$ is a sofic group, it is shown in \cite{israel} that
every injective linear cellular automaton $\tau \colon V^G \to V^G$ is surjective.
The problem of
the existence of an injective but non-surjective linear cellular automaton $\tau \colon V^G \to V^G$, with $V$ a finite-dimensional vector space, remains also open.
\par
If $V$ is a non-zero vector space over an arbitrary field $\K$ and $G$ is a non-periodic group, then there always exists a linear cellular automaton
$\tau \colon V^G \to V^G$ which is surjective but not injective (Proposition \ref{p:lca-inj-not-surj}).
It follows that if $V$ is a non-zero vector 
space over an arbitrary field $\K$ and  $G$ is a group 
such that every surjective linear cellular automaton $\tau \colon V^G \to V^G$ is injective, then $G$ is necessarily periodic (Corollary \ref{c:surj-inj-implies-per}). On the other hand, if $V$ is a finite-dimensional vector space  over an arbitrary field $\K$
and $G$ is a locally finite group, then every surjective linear cellular automaton $\tau \colon V^G \to V^G$ is injective (Proposition \ref{p:loc-fin-surj-lca-implies-inj}). We are able to establish the converse implication  for non-zero finite-dimensional  vector spaces when the ground field is  the field $ \R$ of real numbers  or the field $\C$ of complex numbers, so that we get the following:

\begin{theorem}
\label{t:loc-finite-lca}
Let $G$ be a group. Let $V$ be a non-zero finite-dimensional vector space over  $\K = \R$ or $\C$. 
Then the following conditions are equivalent:
\begin{enumerate}[\rm (a)]
\item
$G$ is locally finite;
\item
every surjective linear cellular automaton $\tau \colon V^G \to V^G$ is injective.
\end{enumerate}
\end{theorem}

As  mentioned above, for an arbitrary ground field $\K$,  
the implication (a) $\Rightarrow$ (b) in Theorem \ref{t:loc-finite-lca}
remains true but we are only able to prove a weak form of the converse implication, namely that (b) implies that $G$ is periodic. 
 This is the reason why we address the following:

\begin{question*}
Let $G$ be a group, $\K$ a field,  and  $V$  a non-zero finite-dimensional vector space over $\K$. 
Suppose that every surjective linear cellular automaton $\tau \colon V^G \to V^G$ is injective.
Does this imply that $G$ is locally finite?
\end{question*}

The proof of Theorem \ref{t:loc-finite-lca}, which is given in Section \ref{sec:proof-main-result}, relies on the techniques of induction and restriction for cellular automata developed in \cite{induction} and on a surjectivity result for the real Laplace operator on finitely generated infinite groups established in \cite{laplace}. Sketches of proofs of these auxiliary results as well as some background material are included in Section \ref{sec:background} for the convenience of the reader.

% SECTION 2
\section{Background material}
\label{sec:background}

The reader is referred to \cite{induction} and \cite{laplace}
for detailed proofs of the statements presented in this section.

\subsection{Restriction of a cellular automaton}
Let $G$ be a group, $A$ a set, and $H$ a subgroup of $G$.
Suppose that a cellular automaton $\tau \colon A^G \to A^G$ admits a memory set $M \subset H$. Let 
$\mu \colon A^M \to A$ denote the associated local defining map.
Then the map $\tau_H \colon A^H \to A^H$ defined by
$$
\tau_H(y)(h) = \mu((h^{-1}y)\vert_M)
\quad \text{ for all } y \in A^H, h \in H,
$$
 is a cellular automaton over the group $H$ and the alphabet $A$ with memory set $M$ and local defining map $\mu$.
 One says that $\tau_H$ is the cellular automaton obtained by  \emph{restriction}  of  
$\tau$ to $H$.
\par
The map $\tau$ may be recovered from $\tau_H$ in the following way.
Let $G/H = \{gH : g \in G\}$ denote the set consisting of all left cosets of $H$ in $G$.
For $c \in G/H$ and $g \in c$, consider the bijective map  $\phi_g \colon H \to c$  defined by $\phi_g(h) = gh$ for all $h \in H$.
Then $\phi_g$ induces a bijective map $\psi_g \colon A^c \to A^H$ given by 
 $\psi_g(z) = z \circ \phi_g$ for all $z \in A^c$.
The cosets $c \in G/H$ form a partition of $G$ so that we can use the identification
$A^G = \prod_{c \in G/H} A^c$.
We then have
\begin{equation}
\label{e;tau-prod}
\tau = \prod_{c \in G/H} \tau_c,
\end{equation}
with $\tau_c \colon A^c \to A^c$ given by $\tau_c = \psi_g^{-1} \circ \tau_H \circ \psi_g$, where $g \in c$.

\begin{proposition}[cf. Th. 1.2 in \cite{induction}]
\label{p:induc-inj-surj}
With the above notation, we have:
\begin{enumerate}[\rm (i)]
\item
$\tau$ is injective if and only if $\tau_H$ is injective;
\item
$\tau$ is surjective if and only if $\tau_H$ is surjective.
\end{enumerate}
\end{proposition}

\begin{proof}
It follows from \eqref{e;tau-prod}
that $\tau$ is injective (resp. surjective) if and only if each $\tau_c$, $c \in G/H$, is injective (resp. surjective).
On the other hand, if we fix $c \in G/H$ and $g \in c$, the fact that the maps $\tau_c$ and $\tau_H$ are conjugate by $\psi_g$  implies that $\tau_c$ is injective (resp. surjective) if and only if $\tau_H$ is injective (resp. surjective). 
\end{proof}

\subsection{The Laplace operator}
Let $G$ be a group and let $S$ be a nonempty finite subset of $G$.
Let $\K = \R$ or $\C$.
 The \emph{Laplace operator} on $G$ with coefficients in $\K$ associated with $S$  
is the   map $\Delta_{G,S}^\K \colon \K^G \to \K^G$ defined by
\begin{equation}
\label{e:laplace}
\Delta_{G,S}^\K(x)(g) = x(g) - \frac{1}{|S|}\sum_{s \in S} x(gs)
\end{equation}
for all $x \in \K^G$ and $g \in G$, where $|S|$ denotes the cardinality of $S$.
It is clear that $\Delta_{G,S}^\K$ is a linear cellular automaton over the group $G$ and the alphabet $\K$ (viewed as a vector space over itself) admitting $S \cup \{1_G\}$ as a memory set. Here we denote by $1_G$ the identity element in $G$.

\begin{theorem}[cf. Th. 1.1 in \cite{laplace}]
\label{t:laplace-surj}
Let $G$ be a  group and let  $S$ be a finite subset of $G$.  
Let $\K = \R$ or $\C$.
Suppose that the subgroup of $G$ generated by $S$ is infinite.
Then the Laplace operator $\Delta_{G,S}^\K \colon \K^G \to \K^G$ is surjective. 
\end{theorem}

\begin{proof}[Sketch of proof]
We can restrict ourselves to the case  $\K = \R$ since $\Delta_{G,S}^\C = \Delta_{G,S}^\R \oplus i\Delta_{G,S}^\R$.
Denoting by $H$ the subgroup of $G$ generated by $S$, it is clear that
$\Delta_{H,S}^\R$ is the restriction of $\Delta_{G,S}^\R$ to $H$.
Thus, by applying Proposition \ref{p:induc-inj-surj}.(ii), we can assume that $S$ generates $G$.
We then distinguish two cases according to whether the group $G$ is amenable or not. 
In the case when $G$ is amenable,
 the surjectivity of $\Delta_{G,S}^\R$ follows from the Garden of Eden theorem for linear cellular automata established in \cite[Th. 1.2]{Eden-linear} since the maximum principle implies that $\Delta_{G,S}^\R$ is pre-injective.
 Suppose now that $G$ is non-amenable.
 Let $\ell^2(G) = \{x \in \R^G: \sum_{g \in G} x(g)^2 <  \infty\}$ denote the Hilbert space of square-summable real-valued functions on $G$.
Then we have $\Delta_{G,S}^{\R}(\ell^2(G)) = \ell^2(G)$ since the Kesten-Day amenability criterion (cf. \cite{kesten}, \cite{day})
implies that $0$ is in the $\ell^2$-spectrum of $\Delta_{G,S}^\R$ in this case.
 On the other hand, it follows from \cite[Lemma 3.1]{Eden-linear} that $\Delta_{G,S}^\R(\R^G)$ is closed in $\R^G$ for the \emph{prodiscrete topology} (that is, the product topology on $\R^G$ obtained by taking the discrete topology on each factor $\R$).
 We conclude by observing that $\ell^2(G)$ is dense in $\R^G$ for the prodiscrete topology.  
\end{proof}

% SECTION 3
\section{Proof of the main result}
\label{sec:proof-main-result}

\begin{proposition}
\label{p:loc-fin-surj-lca-implies-inj}
Let $G$ be a locally finite group and 
let $V$ be a finite-dimensional vector space over an arbitrary field $\K$.
Let $\tau \colon V^G \to V^G$ be a linear cellular automaton.
Then $\tau$ is surjective if and only if  it is injective.
 \end{proposition}

\begin{proof}
Let $M \subset G$ be a memory set for $\tau$.
As $M$ is finite and $G$ is locally finite, the subgroup $H \subset G$  generated by $M$ is finite.
Consider the linear cellular automaton $\tau_H \colon V^H \to V^H$ obtained from $\tau$ by restriction.
 Since $V^H$ is finite-dimensional,   $\tau_H$ is surjective if and only if it is injective.
 On the other hand, it follows from Proposition \ref{p:induc-inj-surj} that $\tau$ is surjective (resp. injective) if and only if $\tau_H$ is surjective (resp. injective).
Thus, $\tau$ is surjective if and only if it is injective.
  \end{proof}

\begin{proof}[Proof of Theorem \ref{t:loc-finite-lca}]
The fact that (a) implies (b) follows from Proposition \ref{p:loc-fin-surj-lca-implies-inj}.
\par
Suppose that $G$ is not locally finite.
Let us show that there exists a linear cellular automaton $\tau \colon V^G \to V^G$ which is surjective but not injective. Let $d = \dim_\K V$.
Without loss of generality, we can assume that $V = \K^d$.
Since $G$ is not locally finite, we can find a finite subset $S \subset G$ such that the subgroup $H$ generated by $S$ is infinite.
Consider the Laplace operator $\Delta_{G,S}^\K \colon \K^G \to \K^G$ defined by
\eqref{e:laplace}.
As mentioned above, $\Delta_{G,S}^\K$ is a linear cellular automaton admitting $S \cup \{1_G\}$ as a memory set. It follows from Theorem \ref{t:laplace-surj} that $\Delta_{G,S}^\K$
  is surjective. Clearly,   $\Delta_{G,S}^\K$ is not injective since all constant configurations are in its kernel.
Consider now the product map $\tau = (\Delta_{G,S}^\K)^d \colon (\K^d)^G \to (\K^d)^G$, where we use the natural identification
$(\K^d)^G = (\K^G)^d$.
Observe that $\tau$ is a linear cellular automaton admitting $S \cup \{1_G\}$ as a memory set.
On the other hand, $\tau$ is surjective but not injective since any product of surjective 
(resp. non-injective) maps is a surjective (resp. non-injective) map.
 This shows that (b) implies (a).
 \end{proof} 

% SECTION 4
\section{Examples of surjective but not-injective linear cellular automata}

\begin{proposition}
\label{p:lca-inj-not-surj}
Let $G$ be a group and let $V$ be a non-zero vector space over a field $\K$.
Suppose that the group $G$ admits an element $g_0$ having infinite order.
Then the map $\tau \colon V^G \to V^G$ defined by
$\tau(x)(g) = x(gg_0)-x(g)$ for all $x \in V^G$ and $g \in G$, is a linear cellular automaton which is surjective but not injective.
 \end{proposition}
 
 \begin{proof}
 It is clear that $\tau$ is a linear cellular automaton admitting $M = \{1_G,g_0\}$ as a memory set. As all constant configurations are mapped to zero, $\tau$ is not injective. 
 \par
 Let us show that $\tau$ is surjective.
  Let $y \in V^G$. Choose a complete set of representatives $ R \subset G$  for the left cosets of the infinite cyclic group generated by $g_0$. This means that every element $g \in G$ can be uniquely written in the form $g = rg_0^n$, where $r \in R$ and $n \in \Z$.
Then,   the configuration $x \in V^G$ defined by
$$
x(rg_0^n)  =
\begin{cases}
0 &\text{ if } n = 0 \\
\sum_{i = 0}^{n - 1}y(rg_0^i)   &\text{ if } n > 0 \\
\sum_{i = 0}^{n - 1}y(rg_0^{n + i})   &\text{ if } n < 0 \\
\end{cases}
$$
for all $r \in R$ and $n \in \Z$,
clearly satisfies $\tau(x) = y$.
Consequently, $\tau$ is surjective.
\end{proof}

\begin{corollary}
\label{c:surj-inj-implies-per}
Let $G$ be a group and let $V$ be a non-zero vector space over a field $\K$.
Suppose that every surjective linear cellular automaton $\tau \colon V^G \to V^G$ is injective.
Then $G$ is periodic.
\qed
\end{corollary}

\begin{remark}
Observe that we can also deduce Proposition \ref{p:lca-inj-not-surj} from Theorem \ref{t:laplace-surj} when $\K = \R$ or $\C$, since $\tau = \Delta_{G,S}^\K$ by taking $S = \{g_0\}$.
In fact, it is easy to check that $\tau$ is pre-injective for any ground field $\K$ so that Proposition \ref{p:lca-inj-not-surj} is a direct consequence of the Garden of Eden theorem for linear cellular automata \cite[Th. 1.2]{Eden-linear}.
 \end{remark}

% REFERENCES

\end{document}